\documentclass[12pt,a4paper]{article}
\pdfoutput=1 %ArXiv
\usepackage[T1]{fontenc}
\usepackage[utf8]{inputenc}
\usepackage[english]{babel}
\usepackage{authblk}
\usepackage{amsmath,amssymb,amsthm}
\usepackage{graphicx}
\usepackage{color}
\usepackage{hyperref}
\usepackage{todonotes}
\usepackage{amsfonts}
\usepackage{mathrsfs}
\usepackage{url}
\newtheorem{theorem}{Theorem}

\newtheorem{lemma}{Lemma}

\begin{document}

\vspace{0.2cm}
\begin{center}
{\bf ON CONDITIONS UNDER WHICH A PROBABILITY\\
DISTRIBUTION IS UNIQUELY DETERMINED BY
ITS MOMENTS}\footnote{The work of E.B. Yarovaya on Sections 2--4 is supported by the Russian Science Foundation (Project no.~19-11-00290) and performed in Steklov Mathematical Institute of Russian Academy of Sciences.}
\end{center}

\vspace{0.4cm}
\begin{center}
{Elena B. Yarovaya\footnote{
%Steklov Mathematical Institute of Russian Academy of Sciences, Moscow, Russia;
Lomonosov Moscow State University, Moscow, Russia; e-mail: yarovaya@mech.math.msu.su},
Jordan M. Stoyanov\footnote{Institute of Mathematics and Informatics, Bulgarian Academy of Sciences, Sofia, Bulgaria; e-mail: stoyanovj@gmail.com},
Konstantin K. Kostyashin\footnote{Lomonosov Moscow State University, Moscow, Russia; e-mail: kostjr@rambler.ru}}
\end{center}

\vspace{0.4cm}\noindent
{\small
 We study the relationship between the well-known Carleman's condition guaranteeing that a probability distribution is uniquely determined by its moments, and a recent easily checkable condition on the rate of growth of the moments. We use asymptotic methods in theory of integrals and involve properties of the Lambert $W$-function to show that the quadratic rate of growth of the ratios of consecutive moments, as a sufficient condition for uniqueness, is more restrictive than Carleman's condition. We derive a series of statements, one of them showing that Carleman's condition does not imply Hardy's condition, although the inverse implication is true. Related topics are also discussed.}

\vspace{0.2cm}\noindent

\vspace{0.2cm}\noindent

\vspace{0.2cm}
\begin{flushright}
\it Dedicated to Albert Nikolaevich Shiryaev\\
on the occasion of his 85th birthday
\end{flushright}

\vspace{0.2cm}
\section{Introduction}

Suppose $X$ is a non-negative random variable with the distribution function $F(x)={\bf P}[X \leq x]$, $x \geq 0$, this is further denoted as: $X \sim F$. We assume that all the moments of the random variable $X$ are finite, that is,
\[
m_{n}= {\bf E}X^n=\int_{0}^{\infty} x^{n}dF(x) < \infty, \qquad n=1,2,\ldots\,.
\]
In the classical setting of the problem of moments the next question emerges: Is the distribution function of the random variable $X$  uniquely determined by the moments $\{m_{n}\}$? In view of a fundamental Kolmogorov's theorem (see~\cite{Shir04:v1:r}) the uniqueness property is attributed to both $X$ and $F$. If the answer is affirmative,  the random variable $X$ is called M-\emph{determinate}  (the terms M-\emph{det},  \emph{determinate} or \emph{moment determinate} \cite{Ahiezer61:e}, \cite{Stoy:e}, \cite{Schmudgen17} are also in use). In such a case we also say that the moment problem has a unique solution. If the answer is negative, that is,  if there are different random variables with the same moment sequence,  $X$ is called M-\emph{indeterminate} (or M-\emph{indet}) and it is said that the moment problem has a non-unique solution. It the case of non-uniqueness there are infinitely many distributions of any kind, all with the same moment sequence (see~\cite{Stoy:JAP04}).

In this paper we do not discuss rather complex criteria of solvability (existence and/or uniqueness) of the moment problem expressed in terms of an infinite sequence of Hankel matrices and the requirement for their positive definiteness (see~\cite{Ahiezer61:e}, \cite{Schmudgen17}, \cite{ShT:r}). For a long time, a special attention has been paid to finding easily checkable conditions even if they are only sufficient or only necessary for either M-determinacy or M-indeterminacy. Since $X \geq 0$, then the support of $F$ is the half-line ${\mathbb R}_+=[0,\infty)$, so  we deal with the Stieltjes moment problem, while for $X$ with  values in the whole line ${\mathbb R}=(-\infty, \infty)$, we deal with the Hamburger moment problem.

It is worth noting that to know if a distribution is M-determinate or M-indeterminate is of interest  by itself.  Moreover, the M-determinacy property is essential in the proof of limit theorems. It is appropriate  to recall the Fr\'echet--Shohat Theorem~\cite{FSh:TAMS31}.
\begin{theorem}%\label{T:FSh}
Let $F_N$, $N=1,2, \ldots$, be a sequence
of distribution functions such that for each of them all moments are finite and the following limits exist:
\[
\lim_{N \to \infty}m_{n,N}=\lim_{N \to \infty}\int x^n\,dF_N(x) = m_n,\qquad n=1,2, \ldots.
\]
Then the following two statements are true:
\begin{enumerate}
\item[\rm(i)] $\{m_n\}$ is a moments sequence of some distribution function, say $F_*$;
\item[\rm(ii)] if $\{m_n\}$  uniquely determines $F_*$, the weak convergence $F_N \stackrel{d}{\to} F_*$, $N \to \infty$, holds.
\end{enumerate}
\end{theorem}

Quite useful is the recent Lin's paper ~\cite{Lin:JSDA17} for providing a systematic description of almost all available conditions, either necessary or sufficient for both M-determinacy and M-indeterminacy. The most important feature of this paper is that the author analyzes the mutual relationships between different conditions. In particular, the most well-known sufficient conditions for M-determinacy, Carleman's condition, is discussed. Let us recall an important theorem in the Stieltjes moment problem (see, e.g.,~\cite[Ch. II, Part `Appendix and Problems']{Ahiezer61:e}, \cite[Section 11, Criterion $(C_{2})$]{Stoy:e}, \cite[Ch. II, \S~12, Section 9]{Shir04:v1:r}).

\begin{theorem}\label{T:Kcond}
Let $\{m_n\}$ be the moment sequence of a random variable $X\ge 0$, and
\begin{equation}\label{E:Carl0}
\sum_{n=1}^\infty {m^{-1/(2n)}_{n}} = \infty\qquad(\textit{Carleman's condition}).
\end{equation}
Then both the distribution function $F$ and the random variable $X$ are uniquely determined by the moments $\{m_{n}\}$,  or, as one says, the random variable $X$ is M-determinate.
\end{theorem}

Though Carleman's condition~\eqref{E:Carl0} (and its analogue $\sum_{n=1}^\infty {m^{-1/(2n)}_{2n}} = \infty$ for distributions on $\mathbb R$) is not necessary for M-determinacy, this condition has been successfully used for characterization of  probability distributions. A series of interesting results can be found, e.g., in the papers~\cite{KlebMk:LOMI81:e}, \cite{Lykov:TVP17:e}, \cite{KLS18,LinSPL:97}. Note also that Carleman's condition may play a key role in the above Fr\'echet-Shohat Theorem when proving limit theorems for branching random walks (see~\cite{HY19:r}, \cite{YarBRW:r}).

Among the variety of available sufficient conditions for  M-determinacy of a random variable $X$ with unbounded support, there are results based on our knowledge of how fast the moments grow, ${\mathbf E}X^n \nearrow \infty$  as $n \to \infty$. In this case it is useful to involve the ratios of two consecutive moments $m_{n+1}/m_n$ (for $X$ with values in ${\mathbb R}_+$) and  $m_{2n+2}/m_{2n}$ (if $X$ takes values in $\mathbb R$), because these ratios characterize the \emph{rate of growth of the moments}. The determinacy of the distribution function $F$ and the random variable $X$ depends essentially on this rate. Recall the following result~\cite[Theorem~1]{LinStoy:JTP15} in the Stieltjes case.

\begin{theorem}\label{T:Stoy}
Let the moments $\{m_k\}$ follow a quadratic rate of growth, that is,
\begin{equation}\label{E:Stoyan0}
\frac{m_{n+1}}{m_{n}}= {\cal{O}}((n+1)^{2}) \quad \text{as} \quad n \to \infty.
\end{equation}
Then the random variable $X$ is M-determinate.
\end{theorem}

A couple of years ago, one of the authors of this paper (J.S.) has given a talk at the Principal Seminar of the Department of Probability Theory (Faculty of Mechanics and Mathematics, Lomonosov Moscow State University). He suggested for M-determinacy to use  the convenient and easily checkable condition~\eqref{E:Stoyan0}  noting that this condition implies  condition \eqref{E:Carl0}, that is,
\[
\eqref{E:Stoyan0} \ \Rightarrow \ \eqref{E:Carl0}.
\]
The question whether the inverse implication $\eqref{E:Carl0} \Rightarrow \eqref{E:Stoyan0}$ is valid was not discussed. As far as we know, this question remained open, and our goal here is to show that in general this implication is not true. We will describe explicitly a positive M-determinate random variable whose moments satisfy Carleman's condition~\eqref{E:Carl0} but not condition~\eqref{E:Stoyan0}. In other words, it will be shown that condition \eqref{E:Stoyan0} is more restrictive than \eqref{E:Carl0}.

In what follows we are using two notations. The first one, $o(1)$ (little $o$) is to denote small numerical quantities (they are different in different places but we do not use specific indices) depending on a real or integer argument and converging to zero as the argument increases to infinity. The second one, ${\mathcal O}(\cdot)$  (big $O$) is to compare two functions $u(t)$, $t>0$ and $v(t)$, $t>0$,  as $t \to \infty$\,: we will write $u(t)=\mathcal{O}(v(t))$, if there is a constant $C \geq 0$ such that $|u(t)| \le C|v(t)|$ for large $t$.

The structure of the paper is as follows. In Section~\ref{S:Main} we formulate one of the main results, Theorem~\ref{T:1}, containing an explicit example of a random variable satisfying Carleman's condition~\eqref{E:Carl0} but not satisfying the growth rate condition~\eqref{E:Stoyan0}. Auxiliary  facts and statements that are necessary to prove Theorem~\ref{T:1} are given in Sections~\ref{S:L1proof} and \ref{S:L2proof}. For the proofs we use techniques from the asymptotic theory of integrals (see \cite{deBruijn:e}, \cite{Fed87:r}), as well as properties of the Lambert $W$-function \cite{CGHJK:96}. In Section~\ref{S:rem},  Theorem~\ref{T:1} is essentially extended, see Theorems~\ref{T:StoyGen} and~\ref{T:5} (Stieltjes case),  Theorem~\ref{T:7} (Hamburger case) and Theorem~\ref{T:8} showing that Hardy's condition is more restrictive than Carleman's condition. At the end, a useful asymptotic for the Euler gamma function is given.

\section{The main result and its proof}\label{S:Main}
One of the goals of this paper is to prove the following theorem.

\begin{theorem}\label{T:1}%Th4
There exists a positive M-determinate random variable $X$  satisfying Carleman's condition~\eqref{E:Carl0} but not satisfying
condition~\eqref{E:Stoyan0}.
\end{theorem}

\vspace{0.1cm}\noindent
{\it Proof.} \ We are using a somewhat modified idea, that was suggested
in~\cite[Section~11.6]{Stoy:e} for comparison of other criteria.

Let $\xi_{1}$, $\eta_{1}$, $\xi_{2}$ and $\eta_{2}$ be independent non-negative random variables each exponentially distributed with parameter $1$, denoted further by Exp. Consider the random variables
\[
X_{1}:=\xi_{1}\ln(1+\eta_{1}), \quad X_{2}:=\xi_{2}\ln(1+\eta_{2}).
\]
Since the variables $\xi_{i}$ and $\eta_{i}$, $i=1,2$,  are independent the expressions for the moments $\mathbf{E}X_{i}^{n}$ can be written as
\[
\mathbf{E} X_{i}^{n}=\mathbf{E} \xi_{i}^{n}\, \mathbf{E} [\ln(1+\eta_{i})]^{n},\qquad n=1,2,\ldots\,.
\]
We use now the fact that the moments of  $\xi_{i} \sim$ Exp are well-known:
\[
\mathbf{E} \xi_{i}^{n}= n!,\qquad n=1,2,\ldots\,.
\]

Now, we define a new random variable,
\begin{equation}\label{E:X}
X:=X_{1}X_{2}
\end{equation}
and will show that all properties, required in Theorem~\ref{T:Kcond}, are fulfilled for this variable.

In view of the independence of the random variables $X_{1}$ and $X_{2}$, we obtain
\begin{equation}\label{E:mn-def}
m_{n}=\mathbf{E} X^{n}=\mathbf{E}X_{1}^{n}\, \mathbf{E}X_{2}^{n} = (n !)^{2}\,K_{n}^{2},
\end{equation}
where
\begin{equation}\label{E:Kn}
K_{n}:=\mathbf{E} [\ln(1+\eta_{1})]^{n}=\int_{0}^{\infty} [\ln(1+x)]^{n}e^{-x}\,dx.
\end{equation}

In order to complete the proof of Theorem~\ref{T:1}, we need one auxiliary notion and two lemmas.

Recall that the Lambert $W$-function, usually denoted by $W(t)$, $t \geq 0$, (see, e.g.,~\cite{CGHJK:96}), is defined as the solution of the  nonlinear functional equation
\[
t=W(t)e^{W(t)},\qquad t \geq 0.
\]
We need a few  properties of the Lambert $W$-function, they are formulated and proved in Section~\ref{S:L1proof}, see Lemma~\ref{L:Lambert}.

\begin{lemma}\label{L:1}
For the moments~\eqref{E:mn-def}, the following equalities are satisfied:
\begin{equation}\label{E:I}
m_{n}= (n!)^{2}\,K_{n}^{2}, \quad n=1,2,\ldots,
\end{equation}
where
\[
K_{n}=e\sqrt{2\pi n}[W(n)]^{n-1/2} e^{-n/W(n)}(1+o(1)), \qquad n\to \infty.
\]
\end{lemma}

\begin{lemma}\label{L:2}
The following relations are true:
\begin{align}\label{E:Kn1n}
K_{n}^{1/n}&=(\ln(n+1))\,(1+o(1)),\qquad n\to\infty,\\
\label{E:KdivK} \frac{K_{n+1}}{K_{n}}&=(\ln(n+1))\,(1+o(1)), \qquad n\to \infty.
\end{align}
\end{lemma}

We postpone the proofs of these lemmas to Sections~\ref{S:L1proof} and~\ref{S:L2proof}, respectively, while now we are in a position to complete the proof of Theorem~\ref{T:1}, assuming that Lemmas~\ref{L:1} and~\ref{L:2} are already proved.

Applying Lemma~\ref{L:1} and relation~\eqref{E:Kn1n} from Lemma~\ref{L:2}, we have
\[
m_{n}^{1/(2n)}=K_{n}^{1/n}\,(n!)^{1/n}=(\ln(n+1))\,(n!)^{1/n} \,(1+o(1)),\qquad n\to\infty.
\]
Here, using the Stirling formula
\[
n! =\sqrt{2\pi n} \left({\frac {n}{e}}\right)^{n}\,(1+o(1)),
\]
we obtain
\begin{equation}\label{E:mnroot}
m_{n}^{1/(2n)}=\frac{n}{e}\,(\ln(n+1))\,(1+o(1)),\qquad n\to\infty,
\end{equation}
which implies Carleman's condition~\eqref{E:Carl0} for the moments~\eqref{E:mn-def}.

From~\eqref{E:I} it follows that
\[
\frac{m_{n+1}}{m_{n}}= (n+1)^{2}\left(\frac{K_{n+1}}{K_{n}}\right)^{2},\qquad n =1,2,\ldots.
\]
This, together with relation~\eqref{E:KdivK} from Lemma~\ref{L:2}, yields
\begin{equation}\label{E:mdivm}
\frac{m_{n+1}}{m_{n}}= (n+1)^2\,(\ln(n+1))^{2}\,(1+o(1)),\qquad n\to\infty,
\end{equation}
and hence condition~\eqref{E:Stoyan0} is not valid for the moments~\eqref{E:mn-def}. Theorem~\ref{T:1} is proved.

\vspace{0.2cm}\noindent
{\bf Remark} \ One of the referees has expressed an opinion that, if we use specific results from sources not in our list of
references, then a shorter proof to Theorem 1 can be provided. We keep our original proof which looks a little long, however
the proof itself is instructive and the intermediate statements are of their own interest.

\section{Proof of Lemma~\ref{L:1}}\label{S:L1proof}
We start with obtaining the necessary properties of the Lambert $W$-function.

\begin{lemma}\label{L:Lambert}%lemma3
For the Lambert $W$-function, the equalities $W(0)=0$ and $W(e)=1$ hold. Moreover, $W(t)>0$ for $t>0$, and the following relations
are satisfied:
\begin{align}\label{E:Wderiv}
W'(t)&=\frac{1}{e^{W(t)}(W(t)+1)}>0, \qquad t\ge 0;\\
\label{E:Lamb-approx0} W(t)&=\ln t -\ln\ln t +\mathcal{O}\left(\frac{\ln\ln t}{\ln t}\right)= \ln t -\ln\ln t +o(1),\qquad t\to\infty.
\end{align}
It follows, in particular, that
\begin{equation}\label{E:Lamb-approx1}
W(t)=(\ln t) (1 +o(1)),\qquad t\to\infty.
\end{equation}
Moreover,
\begin{align}\label{E:WdivW}
\left(\frac{W(t+1)}{W(t)}\right)^{t}&\to 1,\qquad t\to\infty,\\
\label{E:difffrac} \frac{t+1}{W(t+1)}-\frac{t}{W(t)}&\to 0, \qquad t\to\infty.
\end{align}
\end{lemma}

\vspace{0.1cm}\noindent
{\it Proof.} \
By writing the definition of the Lambert $W$-function,
\begin{equation}\label{E:Wdef}
t=W(t)e^{W(t)},
\end{equation}
we easily find that $W(0)=0$ and $W(e)=1$ and also conclude the strict positivity of $W(t), t>0$ and the  inequality~\eqref{E:Wderiv}. Therefore
\begin{equation}\label{E:Wineq}
W(t)>1,\quad \ln W(t)>0\quad\text{for}\quad t>e.
\end{equation}

Taking the logarithm of both sides of~\eqref{E:Wdef}, we get
\begin{equation}\label{E:WlW}
W(t)=\ln t -\ln W(t),
\end{equation}
which implies, in view of~\eqref{E:Wineq}, that
\[
W(t)\le \ln t,\qquad t>e.
\]
Substituting this bound for $W(t)$ in the right side of~\eqref{E:WlW} implies that
\begin{equation}\label{E:Wlow}
W(t)\ge \ln t -\ln\ln t,\qquad t>e.
\end{equation}
A substitution of the bound \eqref{E:Wlow} for $W(t)$ in the right side of~\eqref{E:WlW} leads to
\begin{align}\label{E:Wup}\notag
W(t)&\leq\ln t -\ln \bigl(\ln t -\ln\ln t\bigr)=
\ln t -\ln\ln t -\ln \left(1 -\frac{\ln\ln t}{\ln t}\right)\\
&= \ln t -\ln\ln t +\mathcal{O}\left(\frac{\ln\ln t}{\ln t}\right),\qquad t\to\infty.
\end{align}

By comparing~\eqref{E:Wlow} and~\eqref{E:Wup}, we conclude that relations~\eqref{E:Lamb-approx0} hold, and this implies, in particular, relation~\eqref{E:Lamb-approx1}\footnote{The limit relation~\eqref{E:Lamb-approx0} is known, see, e.g.,~\cite[Section 2.4]{deBruijn:e}. Its proof is given here for the sake of completeness of the exposition.}.

Let us now prove relation~\eqref{E:WdivW}. In view of~\eqref{E:Wderiv}, the function $W(t)$, $t\ge 0$ is increasing, so
\begin{equation}\label{E:WdivWge1}
1\le \frac{W(t+1)}{W(t)}.
\end{equation}

By the mean value theorem, $W(t+1)-W(t)=W'(t+\theta)$, where $\theta$ is a number from the interval $(0,1)$. The increasing of the function $W(t)$, $t \ge 0$, and  equality~\eqref{E:Wderiv} imply that the  derivative $W'(t)$ is decreasing. Hence
\begin{equation}\label{E:WminusW}
W(t+1)-W(t)=W'(t+\theta)\le W'(t)=\frac{1}{e^{W(t)}(W(t)+1)}=\frac{W(t)}{t(W(t)+1)},
\end{equation}
where the last equality follows from~\eqref{E:Wdef}.

From~\eqref{E:WdivWge1} and~\eqref{E:WminusW}, we obtain that
\begin{align*}
1&\le\left(\frac{W(t+1)}{W(t)}\right)^{t}=
\left(1+\frac{W(t+1)-W(t)}{W(t)}\right)^{t}\le\\
&\le \left(1+\frac{1}{t(W(t)+1)}\right)^{t} =\exp\left(t\ln\left(1+\frac{1}{t(W(t)+1)}\right)\right).
\end{align*}
This implies, in view of the easy inequality $\ln(1+x)\le x$ for all $x\ge 0$, that
\begin{align*}
1&\le\left(\frac{W(t+1)}{W(t)}\right)^{t}\le
\exp\left(t\,\ln\left(1+\frac{1}{t(W(t)+1)}\right)\right)\le\\
&\le \exp \left(t\,\frac{1}{t\,(W(t)+1)}\right) = \exp\left(\frac{1}{W(t)+1}\right)\to 1,\qquad t\to\infty.
\end{align*}
These relations prove~\eqref{E:WdivW}.

It remains to establish relation~\eqref{E:difffrac}. From~\eqref{E:Wdef} and the increasing of the function $W(t)$, $t \ge 0$, it follows that
\begin{equation}\label{E:rstbelow}
\frac{t+1}{W(t+1)}-\frac{t}{W(t)}=e^{W(t+1)}-e^{W(t)}\ge 0.
\end{equation}
On the other hand,
\[
\frac{t+1}{W(t+1)}-\frac{t}{W(t)}=\frac{1}{W(t+1)} +t\left(\frac{1}{W(t+1)}-\frac{1}{W(t)}\right).
\]
The last term in this equality is non-positive since the function $W(t)$ is increasing, and that is why
\begin{equation}\label{E:rstabove}
\frac{t+1}{W(t+1)}-\frac{t}{W(t)}\le \frac{1}{W(t+1)} \to 0,\qquad t\to\infty.
\end{equation}

The required relation~\eqref{E:difffrac} follows now from~\eqref{E:rstbelow} and~\eqref{E:rstabove}. Lemma~\ref{L:Lambert} is proved.

Let us move on studying the asymptotic behaviour of $K_{n}$ for large $n$. For this, we first prove an auxiliary lemma about the asymptotic  of $K_{n}$ in a more general setting,  instead of a natural argument $n$, we take the argument to be a positive real number $t\in \mathbb{R}_+$.
\begin{lemma}\label{L:6}
For the function
\begin{equation}\label{E:Stdef}
S(t):=\int_{0}^{\infty} [\ln(1+x)]^{t}e^{-x}\,dx, \qquad t \in {\mathbb R}_+,
\end{equation}
the following relation is true:
\[
S(t)=e\sqrt{2\pi t}\,[W(t)]^{t-1/2}\,e^{-t/W(t)}(1+o(1)), \qquad t\to \infty.
\]
\end{lemma}

\vspace{0.1cm}\noindent
{\it Proof.}
Let us rewrite~\eqref{E:Stdef} in the form
\[
S(t)=\int_{0}^{\infty} e^{Q(x,t)}\,dx,
\]
where
\begin{equation}\label{E:Q}
Q(x,t):=t\ln\ln(1+x)-x.
\end{equation}

In order to find the asymptotic behaviour of the function $S(t)$ as $t \to \infty$, we use the Laplace saddle-point method. We will need the following lemma which combines the statements of Theorem~2.3 and Corollary~2.1 from~\cite{Fed87:r}.
\begin{lemma}\label{L:3}%lemma5
Let $Q(t,x)$, $x >0$, $t> 0$ be a real-valued function for which the following three conditions hold:
\begin{enumerate}
\item[\rm 1)] for all large $t$, with respect to the argument $x \in (0,\infty)$, \ $Q(t,x)$ has a maximum at a point $x_{t}$, which is non-degenerate in the sense that $Q_{xx}''(x_{t},t)\neq 0$, and there exists a function $\mu(t) \to \infty$, $t\to \infty$, such that uniformly in $x$, for
    \[
    x\in U(x_{t}):=\{ x:\, |x-x_{t}| \leq \mu(t)\,|Q_{xx}''(x_{t},t)|^{-1/2} \},
    \]
    we have the relation
\[
Q_{xx}''(x,t)=Q_{xx}''(x_{t},t)(1+o(1)), \qquad t \to \infty;
\]

\item[\rm 2)]\sloppy for each fixed $t>0$, the function $Q(t,x)$ is strictly convex upward, that is  $Q_{xx}''(x,t)<0$ for all  $x>0$ and $t>0$;

\item[\rm 3)]  $\lim_{t\to \infty} x_{t} \sqrt{|Q_{xx}''(x_{t},t)|}=\infty$.
\end{enumerate}
Then, under conditions 1)--3), for the integral
\[
S(t)=\int_{0}^{\infty}e^{Q(x,t)}\,dx, \qquad t>0,
\]
the following asymptotic relation holds:
\begin{equation}\label{E:asymt}
S(t)=\sqrt{-\frac{2\pi}{Q_{xx}''(x_{t},t)}}\,e^{Q(x_{t},t)}(1+o(1)), \qquad t\to \infty.
\end{equation}
\end{lemma}

We continue the proof of Lemma~\ref{L:6} by showing that all conditions in Lemma~\ref{L:3} are satisfied for the function $Q$ defined by~\eqref{E:Q}.

\emph{Verification of Condition \emph{1)} in Lemma~\emph{\ref{L:3}}.} Let us denote by $x_t$ the point, at which the function $Q(x,t)$ reaches the maximum in the argument $x$. To find~$x_t$, we need to solve the equation
\[
Q_{x}'(x,t)=(t\ln\ln(1+x)-x)_{x}'=0.
\]
After differentiating $Q(x,t)$, see~\eqref{E:Q},  we get
\begin{equation}\label{E:t}
t=(1+x)\ln(1+x).
\end{equation}
Therefore the solution $x=x_{t}$ of \eqref{E:t} has the form
\begin{equation}\label{E:xt}
x_{t}=e^{W(t)}-1,
\end{equation}
that is, $x_t$ is expressed via the Lambert $W$-function, defined as the solution of equation \eqref{E:Wdef}. Substituting $x_{t}=e^{W(t)}-1$ in~\eqref{E:Q}, we obtain
\begin{equation}\label{E:Qx}
Q(x_{t},t)=t\ln\ln(1+x_{t})-x_{t}=t\ln W(t)-e^{W(t)}+1.
\end{equation}

Differentiating~\eqref{E:Q} twice with respect to $x$, we find that
\begin{equation}\label{E:Qdoublediff}
Q''_{xx}(x,t)=-\frac{t[1+\ln(1+x)]}{[(1+x)\ln(1+x)]^{2}}.
\end{equation}
Then, in view of~\eqref{E:xt} and~\eqref{E:Qdoublediff},
\begin{equation}\label{E:Qdiffxt}
Q''_{xx}(x_{t},t)=-\frac{t[1+\ln(1+x_{t})]}{[(1+x_{t})\ln(1+x_{t})]^{2}}= -\frac{1+W(t)}{t}.
\end{equation}

Let $x\in U(x_{t})$, that is
\[
|x-x_{t}|\le \mu(t)\,|Q_{xx}''(x_{t},t)|^{-1/2}=\mu(t)\,\sqrt{\frac{t}{1+W(t)}}.
\]
Introduce the function
\[
h(x,t)=\frac{x-x_{t}}{1+x_{t}}.
\]
Taking into account that $1+x_{t}=e^{W(t)}$ and $\frac{t}{W(t)}=e^{W(t)}$, we conclude that
\begin{align*}
|h(x,t)|&\le \mu(t)\,\frac{1}{1+x_{t}}\,\sqrt{\frac{t}{1+W(t)}}\le
\mu(t)\,e^{-W(t)}\,\sqrt{\frac{t}{W(t)}}\\
&=\mu(t)\,e^{-W(t)}\sqrt{e^{W(t)}}=\mu(t)\,e^{-(1/2)W(t)},\qquad x\in U(x_{t}).
\end{align*}

Choose a function $\mu(t)\to\infty$ such that
\[
\mu(t)e^{-(1/2)W(t)}\to 0,\qquad t\to\infty.
\]
If among the many options we take, e.g.,
\begin{equation}\label{E:mudef}
\mu(t)=e^{(1/4)W(t)},
\end{equation}
we obtain
\[
h(x,t)=\mathcal{O}\left(e^{-(1/4)W(t)}\right)=o(1)\to 0,\qquad t\to\infty,
\]
uniformly in $x\in U(x_{t})$.

Now we study the behaviour of $Q''_{xx}(x,t)$ for $x\in U(x_{t})$ and $t\to\infty$. First, we find that
\begin{align*}
Q''_{xx}(x,t)&=-\frac{t\,[1+\ln(1+x)]}{[(1+x)\ln(1+x)]^{2}}\\
&=-\frac{t\,[1+\ln(1+x_{t}+(x-x_{t}))]}{[(1+x_{t}+(x-x_{t}))\ln(1+x_{t}+(x-x_{t}))]^{2}}\\
&=-\frac{t\,\bigl[1+\ln\bigl((1+x_{t})(1+h(x,t))\bigr)\bigr]}
{\bigl[(1+x_{t})(1+h(x,t))\ln\bigl((1+x_{t})(1+h(x,t))\bigr)\bigr]^{2}}\\
&=-\frac{t\,\bigl[1+\ln(1+x_{t})+\ln(1+h(x,t))\bigr]} {\bigl[(1+x_{t})(1+h(x,t))\bigl(\ln(1+x_{t})+\ln(1+h(x,t))\bigr)\bigr]^{2}}.
\end{align*}
As it was mentioned above, $h(x,t)=o(1)$ as $t \to \infty$ uniformly in $x\in U(x_{t})$, and then uniformly in $x\in U(x_{t})$ we have that $\ln(1+p(x,t))=\mathcal{O}(h(x,t))=o(1)$. Hence,
\begin{align}\notag
Q''_{xx}(x,t)&=-\frac{t\left[1+\ln(1+x_{t})+o(1)\right]}
{\bigl[(1+x_{t})(1+o(1))\bigl(\ln(1+x_{t})+o(1)\bigr)\bigr]^{2}}\\
\notag &=-\frac{t\bigl[1+\ln(1+x_{t})]\left[1+\frac{o(1)}{1+\ln(1+x_{t})}\right]}
{\left[(1+x_{t})(\ln(1+x_{t}))(1+o(1))\left(1+\frac{o(1)}{\ln(1+x_{t})}\right)\right]^{2}}\\
\label{E:almostfin} &=-\frac{t[1+\ln(1+x_{t})]}{[(1+x_{t})\ln(1+x_{t})]^{2}}\, \frac{1+\frac{o(1)}{1+\ln(1+x_{t})}} {\left[(1+o(1))\left(1+\frac{o(1)}{\ln(1+x_{t})}\right)\right]^{2}}.
\end{align}
Here, due to~\eqref{E:Qdiffxt},
\[
\frac{t[1+\ln(1+x_{t})]}{[(1+x_{t})\ln(1+x_{t})]^{2}}= -Q''_{xx}(x_{t},t)
\]
and, moreover,
\[
\frac{o(1)}{1+\ln(1+x_{t})}=o(1),\quad \frac{o(1)}{\ln(1+x_{t})}=o(1),\qquad t \to \infty.
\]
Thus~\eqref{E:almostfin} implies the  equalities
\[
Q''_{xx}(x,t)= Q''_{xx}(x_{t},t)\frac{1+o(1)}{\bigl[(1+o(1))(1+o(1))\bigr]^{2}}= Q''_{xx}(x_{t},t)(1+o(1)),\qquad  t\to\infty,
\]
in which, let us recall, all quantities $o(1)$ depend on $x$ and $t$ and tend to zero as $t\to\infty$ uniformly in $x\in U(x_{t})$.

This shows that Condition~1) is fulfilled if we choose the function $\mu(t)$ in the form of~\eqref{E:mudef}. It is easy to see that this is also true for many other choices of  $\mu(t)$.

\emph{Verification of Condition \emph{2)} in Lemma~\emph{\ref{L:3}}.} According to~\eqref{E:Qdoublediff}, we have
\[
Q_{xx}''(x,t)=-\frac{t(1+\ln(1+x))}{(1+x)^{2}\ln^{2}(1+x)}.
\]
Hence, $Q_{xx}''(x,t)<0$ for all $x>0$ and $t >0$, that is Condition~2) in Lemma~\ref{L:3} is satisfied.

\emph{Verification of Condition \emph{3)} in Lemma~\emph{\ref{L:3}}.} By the definition~\eqref{E:Wdef} of $W(t)$ function,
\begin{equation}\label{E:rW}
e^{W(t)}=\frac{t}{W(t)}\to\infty,\qquad t\to\infty.
\end{equation}
Therefore
\[
x_{t}=e^{W(t)}-1=\frac{t}{W(t)}-1=\frac{t}{W(t)}(1+o(1))\to\infty,\qquad t\to\infty,
\]
and then
\begin{align*}
\lim_{t\to \infty} x_{t} \sqrt{|Q_{xx}''(x_{t},t)|}&=
\lim_{t\to \infty}(e^{W(t)}-1)\sqrt{\frac{1+W(t)}{t}}\\
&=\lim_{t \to \infty}\left(\frac{t}{W(t)}-1\right)\sqrt{\frac{1+W(t)}{t}}= \lim_{t \to \infty}\sqrt{\frac{t}{W(t)}}\,(1+o(1))=\infty,
\end{align*}
that is Condition~3) in Lemma~\ref{L:3} is satisfied.
\medskip

Thus, we have shown that all conditions in Lemma~\ref{L:3} are met and then one can use equality~\eqref{E:asymt} from which, in combination with~\eqref{E:Qx} and~\eqref{E:Qdiffxt}, to derive the following:
\begin{align*}
S(t)&= \sqrt{\frac{2\pi t}{W(t)+1}}\,\exp\left(t\ln
W(t)-e^{W(t)}+1\right)(1+o(1))\\
&=e\,\sqrt{\frac{2\pi t}{W(t)+1}}\,[W(t)]^{t}\,\exp\left(-e^{W(t)}\right)(1+o(1)),\qquad t\to \infty.
\end{align*}
Hence, taking into account \eqref{E:rW} and the fact that $W(t)\to \infty$ as $t\to \infty$, the obtained relation can be rewritten in the following form:
\begin{align*}
S(t)&= e\,\sqrt{\frac{2\pi
t}{W(t)}}\,[W(t)]^{t}\,\exp\left(-e^{W(t)}\right)(1+o(1))\\
&=e\,\sqrt{2\pi t}\,[W(t)]^{t-1/2}\,e^{-t/W(t)}(1+o(1)),\qquad t\to \infty.
\end{align*}
Lemma~\ref{L:6} is proved.

Notice, that $K_{n}=S(n)$. Then, referring to Lemma~\ref{L:6}, we arrive at the claim in Lemma~\ref{L:1}.

\section{Proof of Lemma~\ref{L:2}}\label{S:L2proof}\sloppy
First, let us prove equality~\eqref{E:Kn1n}. According to Lemma~\ref{L:1},
\begin{equation}\label{E:K1n}
K_{n}^{1/n}=W(n)\,e^{1/n}(2\pi n)^{1/n}\,[W(n)]^{-1/(2n)}\,e^{-1/W(n)}(1+o(1))^{1/n}.
\end{equation}
Here, the factors $e^{1/n}$, $(2\pi n)^{1/n}$ and $(1+o(1))^{1/n}$ clearly tend to $1$ as $n\to\infty$. Moreover, in view of statement~\eqref{E:Lamb-approx1} from Lemma~\ref{L:Lambert},
\begin{equation}\label{E:Wn}
W(n)=(\ln n)\,(1+o(1)),
\end{equation}
and both factors $[W(n)]^{-1/(2n)}$ and $e^{-1/W(n)}$ also tend to $1$ as $n\to\infty$. Finally, since equality~\eqref{E:Wn} implies that
\[
W(n)=[\ln(n+1)]\, (1 +o(1)),
\]
then the required relation~\eqref{E:Kn1n} follows from~\eqref{E:K1n}.

It remains to establish relation~\eqref{E:KdivK}. By taking $t=n$ in Lemma~\ref{L:6} we find that
\begin{align*}
\frac{K_{n+1}}{K_{n}}&=\frac{S(n+1)}{S(n)}=
\sqrt{\frac{n+1}{n}}\,\frac{[W(n+1)]^{n+1-1/2}}{[W(n)]^{n-1/2}}\\
&\qquad\times \exp\left(\frac{n}{W(t)}-\frac{n+1}{W(n+1)}\right)(1+o(1)).
\end{align*}
Let us rewrite this equality in the form
\begin{align}\label{E:4}\notag
\frac{K_{n+1}}{K_{n}}&=
W(n+1)\sqrt{1+\frac{1}{n}}\left(\frac{W(n+1)}{W(n)}\right)^{n-1/2}\\
&\qquad\times \exp\left(\frac{n}{W(n)}-\frac{n+1}{W(n+1)}\right)(1+o(1)).
\end{align}
There are four factors in \eqref{E:4}. The third and the fourth factors tend to $1$ as $n\to \infty$, which follows from relations \eqref{E:WdivW} and \eqref{E:difffrac}
in Lemma~\ref{L:Lambert}, respectively. The second factor obviously tends to $1$, as well. Hence all factors, except the first one, tend to $1$ and \eqref{E:4} can be rewritten as
\[
\frac{K_{n+1}}{K_{n}}=W(n+1)(1+o(1)), \qquad n\to \infty.
\]
Here, in view of statement~\eqref{E:Lamb-approx1} in Lemma~\ref{L:Lambert}, one has
\[
W(n+1)=\left(\ln(n+1)\right) (1 +o(1)),
\]
which implies the validity of relation~\eqref{E:KdivK}. Lemma~\ref{L:2} is proved.

\section{Further results and remarks}\label{S:rem}

\vspace{0.2cm}\noindent
{\bf Remark 1} \
In proving Theorem~\ref{T:1}, we have used relation~\eqref{E:Kn1n} in Lemma~\ref{L:2} in order to check whether the moments of the random variable $X$ defined by equality~\eqref{E:X} satisfy  Carleman's condition~\eqref{E:Carl0}.  Let us describe a different approach to show the same property,  that the moments of $X$ satisfy Carleman's condition~\eqref{E:Carl0}.

We notice first that the following extension of Theorem~\ref{T:Stoy} is true.
\begin{theorem}\label{T:StoyGen}%Th5
Let for the moment sequence $\{m_{n}\}$ of the random variable $X$ the following relation hold:
\begin{equation}\label{E:StoyanGen}
\frac{m_{n+1}}{m_{n}}= \mathcal{O}\left((n+1)^{2}\,q^2(n+1)\right) \quad \text{as} \quad n\to\infty,
\end{equation}
where $q(n) >0$, $n=1,2,\ldots,$ is a function satisfying the condition
\[
\sum_{n=1}^\infty \frac{1}{n\,q(n)} = \infty.
\]
Then the moments $\{m_n\}$ satisfy Carleman's condition~\eqref{E:Carl0} and by Theorem~\ref{T:Kcond} these moments uniquely determine the random variable $X$.
\end{theorem}
The proof of this theorem almost literally repeats the proof of Theorem~1 from~\cite{LinStoy:JTP15}; therefore, no details are given here. It remains to refer to relations~\eqref{E:mdivm}, that were derived in the proof of Theorem~\ref{T:1}, and see that condition~\eqref{E:StoyanGen} for the moments $m_{n}$ of the random variable $X$ is satisfied with the function $q(n)=\ln n$. Then by Theorem~\ref{T:StoyGen} the random variable $X$  satisfies Carleman's condition~\eqref{E:Carl0}.

\vspace{0.2cm}\noindent
{\bf Remark 2} \
Let us provide the results of some calculations in order to illustrate the behaviour of the ratio
\[
\frac{m_{n+1}}{m_{n}}\,\frac{1}{(n+1)^{2}}=\left(\frac{K_{n+1}}{K_{n}}\right)^{2}
\]
as $n$ increases (see Theorem~\ref{T:Stoy}). Our calculations look as follows:
\[
\frac{K_{1}}{K_{0}}\approx 0.60, \quad \frac{K_{2}}{K_{1}} \approx 0.89, \quad \frac{K_{3}}{K_{2}}\approx 1.09, \quad \frac{K_{4}}{K_{3}}\approx 1.24,\quad \frac{K_{100}}{K_{99}}\approx 3.39,\quad\ldots\,.
\]
They show that $\left(K_{n+1}/K_{n}\right)^{2}$ is increasing with $n$, however the increase is `quite slow'. At the same time the quantities $K_{n}$, $n\geq 2$ themselves grow `rather fast', have a look at specific values:
\begin{align*}
K_{2}&=\int_{0}^{\infty} [\ln(1+x)]^{2}e^{-x}\,dx \approx 0.53,\\
&\ldots\\
K_{99}&= \int_{0}^{\infty} [\ln(1+x)]^{99} e^{-x}\,dx\approx 1.32 \times
10^{41},\\
K_{100}&= \int_{0}^{\infty} [\ln(1+x)]^{100}e^{-x}\,dx\approx 4.47 \times 10^{41}.
\end{align*}
As a corollary, this confirms also the `fast growth' of the moments $m_{n}= (n!\,K_{n})^{2}$ themselves. Here are some calculated values:
\begin{multline*}
m_{1}=1,~ m_{2}=1.13,\ldots,~  m_{99}=1.73 \times 10^{82} \times (99!)^{2},~ m_{100}=2 \times 10^{83} \times (100!)^{2}.
\end{multline*}
The calculations were performed with the help of the package WolframAlpha$^{\circledR}$, whereas for the integrals we used their representations via the Meijer $G$-function~\cite{GR63:e}. It is useful to note that Carleman's condition~\eqref{E:Carl0} is still satisfied even for such a `fast' rate of growth of the moments.

\vspace{0.2cm}\noindent
{\bf Remark 3} \
The arguments used above allow us to formulate a statement which is more general than Theorem~\ref{T:1}. Recall that we are dealing with independent random variables $\xi_1$, $\xi_2$, $\eta_1$, $\eta_2$, all Exp-distributed. Fix real numbers $r_1 \in [0,1]$ and $r_2 \in [0,1]$  and consider the random variables
\[
X_1(r_1):= \xi_1\,(\ln (1+\eta_1))^{r_1},  \quad X_2(r_2):= \xi_2\,(\ln (1+\eta_2))^{r_2}
\]
and their product
\[
X(r_1,r_2):= X_1(r_1)\,X_2(r_2).
\]
Notice that $X(1,1)$ is precisely the random variable $X$ (see equality~\eqref{E:X}), for which Theorem~\ref{T:1} was proved. Using the same notation $m_n$, as before, we find that
\begin{align*}
m_n&={\bf E}[X(r_1,r_2)]^n = (n!)^2\,K_n(r_1)\,K_n(r_2),\\
K_{n}(r_j)&=\mathbf{E} [\ln(1+\eta_{j})]^{nr_j}=\int_{0}^{\infty} [\ln(1+x)]^{nr_j}e^{-x}\,dx, \qquad j=1, 2.
\end{align*}

\begin{theorem}\label{T:5}
For an arbitrary choice of  $r_1 \in [0,1]$ and  $r_2 \in [0,1]$, the moments $\{m_n\}$ of the random variable  $X(r_1,r_2)$ do not satisfy  condition \eqref{E:Stoyan0}, however, they satisfy condition \eqref{E:Carl0}  (hence, all variables $X(r_1,r_2)$ are uniquely determined by their moments).
\end{theorem}

The proof follows the same scheme as in the proof of Theorem~\ref{T:1}. We start with the definition~\eqref{E:Stdef} of the function $S(t)$ implying that $K_n(r_j)=S(nr_j)$. The next reasoning is clear, so details are omitted.

\vspace{0.2cm}\noindent
{\bf Remark 4} \
In Theorem~\ref{T:5}, let us  look at the special case  $r_1=0$, $r_2=1$. We conclude that Theorem~\ref{T:1} is true for the random variable $\xi_1\,\xi_2\,\ln (1+\eta_2)$, which is a `simpler' product than~\eqref{E:X}. Notice, here the logarithmic function as a factor appears only once.

One can take $\xi_1 \stackrel{\text{a.e.}}{=} \xi_2$ and write simply $\xi$ instead of $\xi_1$ and $\xi_2$; also one can take $\eta_1 \stackrel{\text{a.e.}}{=} \eta_2$ and write  $\eta$ instead of $\eta_1$ and $\eta_2$. Then we see that Theorems~\ref{T:1} and~\ref{T:5} will be valid for a wide range of random variables, in fact, for any member of the family
\[
\xi^{\delta}\,(\ln (1+\eta))^r,\qquad r \in [0,1],~ \delta \in [0,2].
\]

Let us note that if $\xi_1, \xi_2$ and $\xi$ are Exp-distributed random variables with independent $\xi_1$ and $\xi_2$, then both conditions~\eqref{E:Carl0} and~\eqref{E:Stoyan0} are satisfies for the product $\xi_1\xi_2$  and for the square $\xi^2$. The difference is that ${\bf E}[\xi_1\xi_2]^n=(n!)^2$, whereas ${\bf E}[\xi^2]^n = (2n)!$.

\vspace{0.2cm}
Note, that the logarithmic function $\ln (1+x)$, $x > 0$ is a slowly varying function, while the power function $x^{\delta}$, $x > 0$ is a regularly varying function. Regarding these classes of functions the reader is referred, e.g., to~\cite{BGT89}. The following natural question arises:

\emph{Is it possible to compare conditions~\eqref{E:Carl0} and~\eqref{E:Stoyan0} for the random variables $X_1:=a(\xi_1)\,b(\eta_1)$ and $X_2:= a(\xi_2)\,b(\eta_2)$, where  $a(x)$, $x>0$ is an arbitrary regularly varying function and $b(x)$, $x>0$ is an arbitrary slowly varying function?}

The material presented in this paper indicates that this can be done, however we do not discuss here this issue.

\vspace{0.2cm}\noindent
{\bf Remark 5} \
Let us consider the random variables $X_1$, $X_2$, $X(r_1,r_2)$ defined above. The analysis provided can be broaden and analogues of Theorems~\ref{T:1} and~\ref{T:5} can be established,  now for distributions on $\mathbb R$  (Hamburger case).

Indeed, let  $Y$ be a random variable taking values in $\mathbb R$ and let ${\mathbf E}|Y|^n<\infty$, $n=1,2,\ldots$. In this case Carleman's condition for the moments $m_n={\mathbf E}Y^n$, $n=1,2,\ldots,$ has the form
\begin{equation}\label{E:KarlemGam}
\sum_{n=1}^{\infty} m_{2n}^{-1/(2n)} = \infty.
\end{equation}
The rate of growth of the moments is now defined not by condition~\eqref{E:Stoyan0}, but in terms of the even order moments (see~\cite{SGD:JSPI14}):
\begin{equation}\label{E:StoyGam}
\frac{m_{2n+2}}{m_{2n}} = {\mathcal O}((n+1)^2), \qquad n \to \infty.
\end{equation}

\begin{theorem}\label{T:7}
One can construct a random variable $Y$ taking values in $\mathbb R$ such that condition~\eqref{E:StoyGam} is not satisfied, however condition~\eqref{E:KarlemGam} is valid.
\end{theorem}

The idea is besides the independent random variables $\xi_1$, $\xi_2$, $\eta_1$, $\eta_2$, each Exp-distributed, to consider an independent from them Bernoulli random variable $\zeta$ with values $+1$ and $-1$, each with probability $1/2$. Set, for example,
\[
Y_1:= \zeta\,\xi_1\,(\ln (1+\eta_1))^{r_1}, \quad Y_2:= \xi_2\,(\ln (1+\eta_2))^{r_2}, \qquad r_1, r_2 \in [0,1].
\]
Then the product  $Y=Y_1\,Y_2$ is a symmetric random variable with values in $\mathbb R$. It is easily seen, that all odd order moments of $Y$ are equal to zero, thus we work only with the even order moments. The following relation turns to be quite useful:
\[
m_{2n}= {\mathbf E}Y^{2n} =(n!)^2\,K_n(r_1)\,K_n(r_2).
\]
Several cases can be considered. For example, if $r_1=1$, $r_2=1$, or $r_1=1$, $r_2=0$, then condition~\eqref{E:KarlemGam} is satisfied, while condition~\eqref{E:StoyGam} is not. If $r_1=0$, $r_2=0$,   both conditions~\eqref{E:KarlemGam} and~\eqref{E:StoyGam} are satisfied. Other possible values of $r_1$ and $r_2$ can be treated in the same way and make respective conclusions.

Finally, let us notice that the random variable $\zeta\,\xi_1$ follows the Laplace distribution (double exponential distribution); its density is $e^{-|x|}/2$, $x \in {\mathbb R}$.

\vspace{0.2cm}\noindent
{\bf Remark 6} \
To continue the previous reasoning, let us compare two more sufficient conditions for M-determinacy.
Recall first the following definition (see~\cite{StoLin:TPA13}). Let $X>0$ be a random variable with all  moments finite. Suppose,
for some $\varepsilon > 0$ the following condition is satisfied:
\begin{equation}\label{E:Hardy}
{\bf E} e^{{\varepsilon}{\sqrt X}} <\infty  \qquad (\emph{Hardy's condition}).
\end{equation}

As shown in the paper~\cite{StoLin:TPA13}, Hardy's condition implies that the random variable $X$  is uniquely
determined by its moments. The question arises of how Hardy's condition~\eqref{E:Hardy} relates to other sufficient conditions
for M-determinacy. In particular, what is the relation between \eqref{E:Hardy} and Carleman's condition? It  should be mentioned here, that Carleman's condition is expressed in terms of  all moments, while Hardy's condition is based on the distribution of $X$.

We have the following result.

\begin{theorem}\label{T:8}
There exists an M-determinate positive random variable $X$ satisfying Carleman's condition~\eqref{E:Carl0}, but not satisfying Hardy's condition~\eqref{E:Hardy}.
\end{theorem}

\vspace{0.1cm}\noindent
{\it Proof.}
The statement in this theorem will be proved for the random variable $X$ introduced by~\eqref{E:X}, which according to Theorem~\ref{T:1}, satisfies Carleman's condition~\eqref{E:Carl0} and therefore is M-determinate,  in view of Theorem~\ref{T:Kcond}.

In the papers~\cite{Lin:JSDA17}, \cite{StoLin:TPA13} it is shown that for any random variable $X \ge 0$, Hardy's condition~\eqref{E:Hardy}
is fulfilled if and only if the moments $\{m_{n}\}$ of $X$ satisfy the inequalities $m_{n}\le (2n)!\,c_{0}^{n}$, $n = 1,2,\ldots$, with some constant $c_{0}$. By Stirling formula, the last inequalities can be rewritten in the following form:
\begin{equation}\label{E:momH}
m_{n}^{1/(2n)}\le [(2n)!]^{1/(2n)}\sqrt{c_{0}}= \frac{2\sqrt{c_{0}}}{e}\,n\,(1+o(1)),\qquad n\to\infty.
\end{equation}
However, it was established in the proof of Theorem~\ref{T:1} that the moments $\{m_{n}\}$ of  $X$ satisfies relation~\eqref{E:mnroot}:
\[
m_{n}^{1/(2n)}=\frac{n}{e}\,(\ln(n+1))\,(1+o(1)),\qquad n\to\infty,
\]
which contradicts~\eqref{E:momH} for any choice of $c_{0}>0$. This contradiction shows that the random variable $X$
does not satisfy Hardy's condition. Theorem~\ref{T:8} is proved.

We would like to emphasize on the fact that in general, the richness of the theory of moment determinacy
is based on the variety of conditions under which a specific property is true or is not true, and
on the relationships between different conditions. This would provide more opportunities for both further development of the theory and the applications. The results derived in this paper can be treated as an addition to the general picture as
described in \cite{KLS18}, \cite{Lin:JSDA17};  see also the references therein. These results, as well as results of other cited authors, confirm that, in a sense, Carleman's condition is `the best' sufficient condition for
a probability distribution to be uniquely determined  by its moments. Additional and useful material can be found in  the following (yet) unpublished work:

{\sc{S. Sodin}}. {\it Lecture Notes on the classical moment problem.} School of Mathematical

Sciences, Queen Mary London University, March 2019. 47 pp.

\vspace{0.2cm}\noindent
{\bf Remark 7.}
It was shown above that the rate of growth of the moments of the random variable $X$ depends on the asymptotic properties of
the numerical sequence $K_n$, $n=1,2,\ldots$ (see formula (\ref{E:Kn})). This is closely related to properties of the Euler gamma function
\begin{equation}\label{E:Gamma}
\Gamma(z)=\int_{0}^{\infty}t^{z-1}e^{-t}\,dt,\qquad z\in\mathbf{C},~\mathop{\mathrm{Re}}z>0.
\end{equation}
The asymptotics of the function $\Gamma(\cdot)$ and of its derivatives are important in various applications.
That is why we are giving below some additional facts which seem to be useful not only for understanding the results in this paper,
but also for applied problems.

In~\cite[Problem 6.48(b)]{BenOrs78} the question is posed to find the main term of asymptotic of the quantity
\[
\gamma_{n}:=\frac{d^{n}}{dz^{n}}\Gamma(z)\big|_{z=1}\quad\text{as}\quad n\to\infty.
\]
To answer this we need  $n$-times differentiation of~\eqref{E:Gamma} in $z$,
\[
\frac{d^{n}}{dz^{n}}\Gamma(z)=\int_{0}^{\infty}(\ln t)^{n}t^{z-1}e^{-t}\,dt,
\]
which implies that
\[
\gamma_{n}=\int_{0}^{\infty}(\ln t)^{n}e^{-t}\,dt.
\]
Represent the integral in the right side as the sum:
\[
\gamma_{n}=\int_{0}^{1}(\ln t)^{n}e^{-t}\,dt+ \int_{1}^{\infty}(\ln t)^{n}e^{-t}\,dt.
\]
Then, substituting $t=1+x$ in the second integral, we get
\begin{equation}\label{E:gamman}
\gamma_{n}=\int_{0}^{1}(\ln t)^{n}e^{-t}\,dt+ e^{-1}\int_{0}^{\infty}[\ln(1+x)]^{n}e^{-x}\,dx.
\end{equation}
Here, the first integral can be estimated as follows:
\[
e^{-1}n!=e^{-1}\left|\int_{0}^{1}(\ln t)^{n}\,dt \right|\le \left|\int_{0}^{1}(\ln t)^{n}e^{-t}\,dt \right|\le \left|\int_{0}^{1}(\ln t)^{n}\,dt \right| = n!.
\]
The second one, due to~\eqref{E:Kn}, can be expressed in the form:
\[
e^{-1}\int_{0}^{\infty}[\ln(1+x)]^{n}e^{-x}\,dx = e^{-1}K_n.
\]

Thus, the first term in~\eqref{E:gamman} is increasing at a rate  $n!$, or, by the Stirling formula,
the rate is  $\sqrt{2\pi n} \left(n/e\right)^{n}$. Because of~\eqref{E:Kn1n}, the second term  increases not faster
than $\left(C\ln n\right)^{n}$. Therefore,
\[
\gamma_{n}=\left(\int_{0}^{1}(\ln t)^{n}e^{-t}\,dt\right)\,(1+o(1)),
\]
meaning that for large $n$ the quantity $\gamma_{n}$ increases at a rate $n!$.

\end{document}